\documentclass[12pt]{amsart}

\usepackage{amssymb}

\newtheorem{theorem}{Theorem}[section]
\newtheorem{lemma}[theorem]{Lemma}

\numberwithin{equation}{section}

\begin{document}

\baselineskip=16.5pt

\title[Parabolic Higgs bundle and Hilbert scheme]{Symplectic structures
on moduli spaces of parabolic Higgs bundles and Hilbert scheme}

\author[I. Biswas]{Indranil Biswas}

\address{School of Mathematics, Tata Institute of Fundamental
Research, Homi Bhabha Road, Bombay 400005, India}

\email{indranil@math.tifr.res.in}

\author[A. Mukherjee]{Avijit Mukherjee}

\address{Max-Planck-Institute for Mathematics in the Sciences,
Inselstr. 22-26, D-04103 Leipzig, Germany}

\email{avijit@mis.mpg.de}

\date{}

\begin{abstract}

Parabolic triples of the form $(E_*,\theta,\sigma)$ are considered,
where $(E_*,\theta)$ is a parabolic Higgs bundle on a given
compact Riemann surface $X$ with parabolic structure on a fixed
divisor $S$, and $\sigma$ is a nonzero section of the underlying
vector bundle. Sending such a triple to the Higgs bundle
$(E_*,\theta)$ a map from the moduli space of stable
parabolic triples to the moduli space of stable parabolic Higgs
bundles is obtained. The pull back, by this map, of the symplectic
form on the moduli space of stable parabolic Higgs bundles will be
denoted by $\text{d}\Omega'$. On the other hand, there
is a map from the moduli space of stable parabolic triples
to a Hilbert scheme
$\text{Hilb}^\delta(Z)$, where $Z$ denotes the total space of the
line bundle $K_X\otimes{\mathcal O}_X(S)$, that sends a triple
$(E_*,\theta,\sigma)$ to the divisor defined by the section $\sigma$
on the spectral curve corresponding to the parabolic Higgs bundle
$(E_*,\theta)$. Using this map and a meromorphic one--form
on $\text{Hilb}^\delta(Z)$, a natural two--form on the
moduli space of stable parabolic triples is constructed. It is
shown here that this form coincides with the above mentioned
form $\text{d}\Omega'$.

\end{abstract}

\maketitle

\section{Introduction} 

In \cite{BM}, we proved that the two--form on the moduli
space of triples of the form $(E,\theta,\sigma)$, where $(E,\theta)$
is a stable Higgs bundle on Riemann surface $X$ and
$\sigma$ a nonzero
section of $E$, obtained by pulling back the natural symplectic
form on the moduli space of stable Higgs bundles coincides
with the pullback of the symplectic structure on the Hilbert
scheme of zero dimensional subschemes of the total space of
of the holomorphic cotangent bundle $K_X$ of $X$.
Our aim here is to establish an analogous result in the context
of parabolic triples.

Let $X$ be a compact connected Riemann surface of genus $g$.
Fix a finite subset $S$ of $X$. A parabolic vector bundle of
rank two with parabolic structure over $S$ consists of a
holomorphic vector bundle $E$, a line $F_s\, \subset\, E_s$,
and a real number $0< \lambda_s< 1$ for each $s\in S$. A Higgs
field on this parabolic bundle is a holomorphic section $\theta$
of the vector bundle
$\text{End}(E)\otimes K_X\otimes{\mathcal O}_X(S)$ such that
$\theta(s)$ is nilpotent with respect to the flag
$0\, =\, F^0 \, \subset\, F^1 \, =\, F_s\, \subset\, F^2 \, =\,
E_s$ of the fiber $E_s$ for each $s\in S$.

We fix the parabolic weights $\{\lambda_s\}_{s\in S}$ and consider
parabolic Higgs bundles of rank two and fixed positive
degree $d$, with
$d\, >\, 6(g-1) +\# S$. Let ${\mathcal M}^s_H$
denote the moduli space of stable parabolic Higgs bundles. For
such stable parabolic Higgs bundles the dimension of the space
of all holomorphic sections of the underlying vector bundle
remains constant over the moduli space (Lemma
\ref{lemma-vanish}). Therefore, there is a projective bundle 
\begin{equation*}
\phi\, :\, {\mathbb P}_H \, \longrightarrow\, {\mathcal M}^s_H
\end{equation*}
whose fiber over any point representing a parabolic
Higgs bundle $((E, \{F_s\}), \theta)$ is the projective
space ${\mathbb P}H^0(X,\, E)$. In other words,
${\mathbb P}_H$ is the moduli space of triples of the
form $((E, \{F_s\}), \theta,\sigma)$, where $((E, \{F_s\}), \theta)$
is a stable parabolic Higgs bundle and $\sigma$ is a nonzero section
of $E$. See \cite{BG} and \cite{Ga}, where the notion of triples were
introduced, for a detailed study and many
interesting results on triples.

The moduli space ${\mathcal M}^s_H$ has a natural holomorphic
one--form, which we call $\Omega$, such that $\text{d}\Omega$
is a symplectic form on ${\mathcal M}^s_H$. The total
space of the holomorphic cotangent bundle of the moduli space
of stable parabolic bundles sits inside ${\mathcal M}^s_H$
as a Zariski open dense subset. The restriction of $\Omega$
to this open subset coincides with the canonical one--form on the
total space of any cotangent bundle. The main result
proved here relates the form $\phi^*\Omega$ on ${\mathbb P}_H$
with a certain one--form on a Hilbert scheme.

Let $Z$ denote the surface defined by the total space of
the line bundle $K_X\otimes{\mathcal O}_X(S)$ over $X$.
Given a stable parabolic
Higgs bundle $((E, \{F_s\}), \theta)$, there is a spectral curve,
which is a divisor on $Z$ and a rank one torsionfree
sheaf $\mathcal L$ on  it such that $\gamma_*{\mathcal L}\,
\cong\, E$, where $\gamma$ is the projection map of the spectral
curve to $X$. Since $\gamma$ is a finite map, we have
$H^0(X,E)$ identified with the space of
all sections of $\mathcal L$. Considering the divisor
of the section of $\mathcal L$ corresponding to $\sigma \in
H^0(X,E)\setminus\{0\}$ we get a map from ${\mathbb P}_H$
to the Hilbert scheme $\text{Hilb}^\delta(Z)$
of zero dimensional subschemes of $Z$ of length $\delta =
d+2(g-1)+\# S$. This map, which we will denote by
$f$, sends a parabolic triple $((E, \{F_s\}), \theta,\sigma)$ to
the divisor on the spectral curve for $((E, \{F_s\}), \theta)$
defined by $\sigma$. Note that using the inclusion map of a spectral
curve in $Z$, a divisor on a spectral curve is a zero
dimensional subscheme of $Z$. (See Section \ref{para-spec.-da}
for the details.)

Using the fact that
$Z$ is the total space of $K_X\otimes{\mathcal O}_X(S)$, there
is a natural meromorphic one--form $\Omega_\delta$ on
$\text{Hilb}^\delta(Z)$. The pullback
$f^*\Omega_\delta$ is a holomorphic one--form on
${\mathbb P}_H$.

We prove that $f^*\Omega_\delta$ coincides with $\phi^*\Omega$
(Theorem \ref{main-theroem}).

Although we have restricted ourselves to rank two parabolic
bundles, the extension of Theorem \ref{main-theroem} to
higher rank case is straight forward. The reason for restriction
to rank two case is the ensuing notational simplification.

\section{Preliminaries}

\subsection{Parabolic Higgs bundles}\label{para.-def.}
Let $X$ be a compact connected Riemann surface of genus
$g$. Fix a finite subset
\begin{equation*}
S \, :=\, \{s_1, s_2, \cdots , s_n\}\, \subset\, X\, .
\end{equation*}
If $g=0$, then take $n\geq 4$, and $n\geq 1$ if $g=1$.

A \textit{parabolic vector bundle of rank two} over $X$
with parabolic structure over $S$ consists
of the following \cite{MS}:
\begin{enumerate}
\item{} a holomorphic vector bundle $E$ of rank two over $X$;
\item{} for each point $s\in S$, a line $F_s\,\subset\, E_s$
of the fiber $E_s$;
\item{} for each point $s\in S$, a real number $\lambda_s$
with $0\,< \, \lambda_s \, <\, 1$.
\end{enumerate} 
The numbers $\{\lambda_s\}$ are called \textit{parabolic weights}.
See \cite{MS} for the definition of parabolic (semi)stability.

\textit{We fix real numbers $\{\lambda_s\}_{s\in S}$ and an
integer $d$. Henceforth, by a parabolic vector bundle we will
always mean a parabolic vector bundle over $X$
of rank two and degree
$d$ with parabolic structure over $S$ and with
parabolic weight $\lambda_s$ for each $s\in S$.}

Let $K_X$ denote the holomorphic cotangent bundle of $X$.
A \text{Higgs structure} on a parabolic vector bundle
$E_*\, :=\, (E, \{F_s\})$ is a holomorphic section
\begin{equation}\label{def.-Higgs-field}
\theta\, \in\, H^0(X,\, \text{End}(E)\otimes K_X\otimes
{\mathcal O}_X(S))
\end{equation}
with the property that for each $s\in S$, the image of
the homomorphism
\begin{equation*}
\theta(s) \, :\, E_s\, \longrightarrow\,
(E\otimes K_X\otimes{\mathcal O}_X(S))_s
\end{equation*}
is contained in the subspace
$F_s\otimes (K_X\otimes{\mathcal O}_X(S))_s\, \subset\,
(E\otimes K_X\otimes{\mathcal O}_X(S))_s$ and
$\theta(s)(F_s) \, =\, 0$ \cite{Hi1}, \cite{BR}, \cite{Fa}.
In other words, $\theta(s)$ is nilpotent with respect
to the flag $0\subset F_s \subset E_s$.

A parabolic Higgs bundle $(E_*\, ,\theta)$ as above is called
\textit{stable} if for every line subbundle $L$ of $E$ with
\begin{equation*}
\theta(L) \, \subseteq\, L\otimes K_X\otimes
{\mathcal O}_X(S) \, \subset\,
E\otimes K_X\otimes{\mathcal O}_X(S)
\end{equation*}
the following inequality is satisfied
$$
\text{degree}(L) + \sum_{s\in S'}\lambda_s \, <\,
\frac{\text{par-deg}(E_*)}{2}\, ,
$$
where $S' \, :=\, \{s\in S\, \vert\, L_s \, =\, F_s\}$; if the
strict inequality is replaced by partial inequality, then
$(E_*\, ,\theta)$ is called \textit{semistable}.

The moduli space of semistable parabolic Higgs bundles
will be denoted by ${\mathcal M}_H$. It is an irreducible
normal quasiprojective variety, and ${\mathcal M}^s_H$ is
a Zariski open smooth subvariety of it.

\subsection{Symplectic structure on moduli
space}\label{subsec.-symp.-str.}
The moduli space ${\mathcal M}^s_H$ has a natural
holomorphic symplectic structure \cite[Section 6]{BR};
we will briefly recall it here.

Take a stable parabolic Higgs bundle $(E_*\, ,
\theta)\, :=\, ((E, \{F_s\}),
\theta)$ represented by a point in ${\mathcal M}^s_H$. Define
\begin{equation*}
\text{End}^1(E) \, \subset\, \text{End}(E)
\end{equation*}
by the condition
that for any $s\in S$ and $v\in \text{End}^1(E)_s$ we have
$v(F_s) \, \subseteq\, F_s$. Let
\begin{equation*}
\text{End}^0(E) \, \subset\, \text{End}^1(E)
\end{equation*}
be defined by the condition $v(F_s)\, =\, 0$.
Consider the two term complex $C_.$ of sheaves
defined by
\begin{equation}\label{2-term-complex}
C_0\, :=\, \text{End}^1(E) \, \stackrel{[-,\theta]}{\longrightarrow}
\, C_1 \, :=\, \text{End}^0(E)\otimes K_X\otimes{\mathcal O}_X(S)
\end{equation}
where $C_i$ is at the $i$th position.

The tangent space $T_{(E_*,\theta)}{\mathcal M}^s_H$ of the variety
${\mathcal M}^s_H$ at the point represented by the parabolic
Higgs bundle $(E_*\, ,\theta)$ is identified with the hypercohomology
${\mathbb H}^1(C_.)$ \cite[Lemma 6.1]{BR}.
Consider the homomorphism
\begin{equation*}
\begin{array}{ccc} C_0 &
\stackrel{[-,\theta]}{\longrightarrow} & C_1\\
\Big\downarrow & & \Big\downarrow\\
\text{End}^1(E) & {\longrightarrow} & 0
\end{array}
\end{equation*}
This induces a homomorphism
\begin{equation}\label{inf.-deform.-forget}
f\, :\, H^1(C_.) \, \longrightarrow\, H^1(X,\, \text{End}^1(E))\, .
\end{equation}
Now observe that $\text{End}^1(E)^*\, \cong\,
\text{End}^0(E)\otimes{\mathcal O}_X(S)$
with the duality pairing defined by
\begin{equation*}
(\omega\, ,\alpha) \, \longmapsto\, \text{trace}(\omega\alpha)
\, \in\, {\mathbb C}\, ,
\end{equation*}
where $\alpha\in \text{End}^1(E)_x$ and
$\omega\in (\text{End}^0(E)\otimes{\mathcal O}_X(S))_x$, and
$x$ is any point of $X$. Now the Serre duality says
\begin{equation*}
H^1(X,\, \text{End}^1(E))^* \, \cong\,
H^0(X,\, \text{End}^0(E)\otimes K_X\otimes
{\mathcal O}_X(S))\, .
\end{equation*}
Consequently, we have a functional $\theta'
\in\, H^1(C_.)^*$ defined by $\beta \, \longmapsto \, (\theta\,
,f(\beta))$, where $f$ is constructed in \eqref{inf.-deform.-forget}.

Let $\Omega$ denote the one--form on the variety ${\mathcal M}^s_H$
that sends any tangent vector $\beta \in T_{(E_*,\theta)}
{\mathcal M}^s_H$, where $(E_*\, ,\theta)\,\in\, {\mathcal M}^s_H$,
to $\theta'(\beta)$ constructed above.

The two--form ${\rm d}\Omega$ is a symplectic form on
${\mathcal M}^s_H$. The restriction of ${\rm d}\Omega$ to the Zariski
open subset of ${\mathcal M}^s_H$ defined by the total space of the
cotangent bundle of the moduli space of parabolic bundles coincides
with the canonical symplectic structure on cotangent bundles.
(See \cite{BR}, \cite{Fa}.)

\subsection{Spectral data for parabolic Higgs
bundles}\label{sec.-spectral-data}
Let $(E_*\, , \theta)$ be a parabolic Higgs bundle. So we have
$\text{trace}(\theta)\, \in\, H^0(X, K_X)$
and
\begin{equation*}
\text{trace}(\theta^2)\, \in\, H^0(X, K^{\otimes 2}_X\otimes
{\mathcal O}_X(S))
\end{equation*}
since $\theta (s)$ is nilpotent for each $s\in S$. Set
\begin{equation}\label{def.-hitchin-space}
{\mathcal H} \, :=\, H^0(X, K_X)\oplus H^0(X, K^{\otimes 2}_X
\otimes{\mathcal O}_X(S))\, .
\end{equation}
Hitchin defined a map
\begin{equation}\label{def.-hitchin-map}
\psi\, :\, {\mathcal M}_H \, \longrightarrow\, {\mathcal H}
\end{equation}
that sends any semistable parabolic Higgs bundle $(E_*\, , \theta)$
to $(\text{trace}(\theta)\, , \text{trace}(\theta^2))\, \in\,
{\mathcal H}$ \cite{Hi2}, \cite{Hi1}, \cite{Fa}. This map
$\psi$ is known as the \textit{Hitchin map}, and
$\mathcal H$ is known as the \textit{Hitchin space}.

Let $Z$ denote the total
space of the line bundle $K_X\otimes{\mathcal O}_X(S)$, which 
is a quasiprojective complex surface. The natural projection of
$Z$ to $X$ will be denoted by $\gamma$. Since $S$ is an effective
divisor, there is a natural homomorphism from $K_X^{\otimes 2}
\otimes{\mathcal O}_X(S)$ to $K^{\otimes 2}_X\otimes
{\mathcal O}_X(2S)$; this homomorphism will be denoted by $q$.

Take any point $(\alpha\, ,\beta)\, \in\,
{\mathcal H}$. Consider the map from $Z$ to the total space of the
line bundle $K^{\otimes 2}_X\otimes{\mathcal O}_X(2S)$ defined by
\begin{equation}\label{spec.-curve-loc.}
z \, \longmapsto\, z\otimes z + q(z\otimes
\alpha (\gamma(z)))+q(\beta (\gamma(z)))\,\in\,
(K^{\otimes 2}_X\otimes{\mathcal O}_X(2S))_{\gamma(z)}\, .
\end{equation}
The inverse image (in $Z$) by this map of the zero section
of $K^{\otimes 2}_X\otimes{\mathcal O}_X(2S)$ is the
\textit{spectral curve} associated to the point
$(\alpha\, ,\beta)$ of the Hitchin space (see \cite{Hi2}, \cite{Fa}).

For any point $h \, =\, (\alpha\, ,\beta)\, \in\,
{\mathcal H}$, the corresponding spectral curve will be denoted
by $Y_h$. The restriction to $Y_h$ of the projection $\gamma$
to $X$ will also be denoted by $\gamma$. Note that the map
\begin{equation}\label{spec.-proj.}
\gamma\, :\, Y_h\, \longrightarrow \, X
\end{equation}
is of degree two. This map is
evidently ramified over every point of $S$.

Given a parabolic Higgs bundle, there is a corresponding
spectral curve and a torsionfree sheaf on it \cite{Hi2}, \cite{Fa}.
To describe the construction, take a semistable parabolic Higgs
bundle $(E_*\, ,\theta)\, \in\, {\mathcal M}_H$
on $X$. The spectral curve associated to
it is the one defined by $\psi((E_*\, ,\theta))\, \in\,
{\mathcal H}$ by the above construction, where $\psi$ is the
Hitchin map defined in \eqref{def.-hitchin-map}. Denote
the point $\psi((E_*\, ,\theta))$ by $h$.

There is a torsionfree sheaf ${\mathcal L}$
of rank one on the spectral curve $Y_h$
such that
\begin{equation}\label{direct-image-l}
\gamma_*{\mathcal L} \, \cong\, E
\end{equation}
(the underlying vector bundle of the parabolic bundle), where
$\gamma$ as in \eqref{spec.-proj.}. The spectral curve
can be thought of as the eigenvalues of the endomorphism
$\theta$. The sheaf ${\mathcal L}$ is defined by the
corresponding eigenvectors (see \cite{Hi2} for the details).
Since the map $\gamma$ is ramified over any point of $s\in S$,
the direct image $\gamma_*{\mathcal L}$ has a filtration over
$s$. This filtration is defined by the order of the vanishing
(at $\gamma^{-1}(s))$ of a (locally defined) section of
$\mathcal L$. In the isomorphism of $\gamma_*{\mathcal L}$
with $E$, The filtration of $\gamma_*{\mathcal L}$ at any
$s\in S$ coincides with the filtration $F_s\,\subset\, E_s$
for the parabolic structure.

\section{Parabolic triples and Hilbert scheme}

In Section \ref{para.-def.} we fixed the degree of a parabolic
vector bundle to be $d$.

\textit{Henceforth, we will assume that the integer $d$, which is
the degree of a parabolic vector bundle, satisfies the condition
\begin{equation*}
d\, > \, 6g-6+n
\end{equation*}
where $n\, =\, \# S$.}

\begin{lemma}\label{lemma-vanish}
For any semistable parabolic Higgs bundle $(E_*\, ,\theta)\, \in\,
{\mathcal M}_H$ over $X$, we have $H^1(X,\, E)\, =\, 0$, where $E$ is
the underlying vector bundle. Consequently, $\dim H^0(X,\, E)\, =\,
d+2(1-g)$.
\end{lemma}

\begin{proof}
Take any $(E_*\, ,\theta)\, \in\, {\mathcal M}_H$. Let $E$ be the
underlying vector bundle for the parabolic vector bundle $E_*$.
Since $H^1(X,\, E)\, \cong\, H^0(X,\, E^*\otimes K_X)^*$, it
suffices to show $H^0(X,\, E^*\otimes K_X)\, =\, 0$. Assume that
$\tau\in H^0(X,\, E^*\otimes K_X)\setminus \{0\}$ is
a nonzero section.

So $\tau$ defines a nonzero homomorphism from $E$ to $K_X$, which will
be denoted by $\overline{\tau}$. The kernel of $\overline{\tau}$ is a
torsionfree coherent subsheaf of $E$. Therefore, it defines a line
bundle over $X$, which will be denoted by $L$. Now we have a diagram
\begin{equation}\label{lemm.-diag}
\begin{array}{ccccc}
L & \stackrel{\iota}{\longrightarrow} & E
& \stackrel{\overline{\tau}}{\longrightarrow} & K_X\\
& & \Big\downarrow\theta && \\
L\otimes K_X\otimes{\mathcal O}_X(S)
& \stackrel{\iota\otimes\text{Id}}{\longrightarrow} & E
\otimes K_X\otimes{\mathcal O}_X(S)
& \stackrel{\overline{\tau}\otimes\text{Id}}{\longrightarrow} &
K^{\otimes 2}_X\otimes{\mathcal O}_X(S)
\end{array}
\end{equation}
We will show that the composition $(\overline{\tau}\otimes\text{Id})
\circ\theta\circ \iota\, =\, 0$. To prove this, first note that
the top sequence in \eqref{lemm.-diag} shows that
\begin{equation}\label{inequality}
\text{degree}(L)\, \geq \, \text{degree}(E)-\text{degree}(K_X)
\, =\, d-2g+2\, .
\end{equation}
On the other hand, $\text{degree}(K^{\otimes
2}_X\otimes{\mathcal O}_X(S))\, =\, 2g-4+n$,
where $n$ is the cardinality of the set $S$. Since $d$ is assumed
to be at least $6g-5+n$, we have
\begin{equation*}
\text{degree}(L)\, > \, \text{degree}(K^{\otimes
2}_X\otimes{\mathcal O}_X(S))\, .
\end{equation*}
Consequently, there is no nonzero homomorphism from $L$ to
$K^{\otimes 2}_X\otimes{\mathcal O}_X(S)$. In particular,
the composition $(\overline{\tau}\otimes\text{Id})
\circ\theta\circ \iota\, =\, 0$.

Let $L'$ denote the line subbundle of $E$ generated by $L$ (note that
$\iota$ may not be fiberwise injective). Since
$(\overline{\tau}\otimes\text{Id})
\circ\theta\circ \iota\, =\, 0$, it follows immediately, that
$\theta(L') \, \subseteq\, L'\otimes K_X\otimes{\mathcal O}_X(S)$.
Finally, we have
\begin{equation*}
\text{degree}(L')\, \geq \, \text{degree}(L) \geq  d-2g+2 =
\frac{d}{2} + \frac{1}{2}(d-4g+4) > \frac{d}{2}+ \frac{n}{2} =
\frac{\text{par-deg}(E_*)}{2}
\end{equation*}
where the second inequality was obtained in \eqref{inequality}
and the third one follows from the assumption
$d\, > \, 6g-6+n$; the first inequality follows from the fact
that there is a nonzero homomorphism from $L$ to $L'$.
The above inequality shows that the line subbundle $L'$
of $E$ contradicts the semistability condition of the parabolic
Higgs bundle $(E_*\, ,\theta)$.

Therefore, we have $H^1(X,\, E)\, =\, 0$. Now the
Riemann--Roch says that $\dim H^0(X,\, E)\, =\,
d+2(1-g)$. This completes the proof of the lemma.
\end{proof}

The above lemma says that $\dim H^0(X,\, E)$ remains constant
over ${\mathcal M}^s_H$. Therefore, there is a natural projective
bundle
\begin{equation}\label{moduli-triple}
\phi\, :\, {\mathbb P}_H \, \longrightarrow\, {\mathcal M}^s_H
\end{equation}
of relative dimension $d-2g+1$ such that the fiber over any
point $(E_*\, ,\theta)\, \in\, {\mathcal M}_H$ is the projective
space ${\mathbb P}H^0(X,\, E)$ consisting of all lines
in $H^0(X,\, E)$; as before, $E$ denotes the underlying vector
bundle for the parabolic bundle $E_*$.

Therefore, ${\mathbb P}_H$ is the moduli space of triples of
the form $(E_*\, ,\theta\, ,\sigma)$, where $(E_*\, ,\theta)$ is a
stable parabolic Higgs bundle and $\sigma\, \in\, H^0(X,\, E)
\setminus \{0\}$ a nonzero section.

By a \textit{parabolic triple} we will mean a triple
$(E_*\, ,\theta\, ,\sigma)$ of the above type. Consequently,
${\mathbb P}_H$ is the moduli space of all parabolic triples.

In Section \ref{subsec.-symp.-str.} we defined a symplectic
structure $d\Omega$ on ${\mathcal M}^s_H$. Define the
algebraic one--form
\begin{equation}\label{triple-one-form}
\Omega'\, :=\, \phi^*\Omega
\end{equation}
on ${\mathbb P}_H$, where $\phi$ is defined in \eqref{moduli-triple}.
So,
\begin{equation*}
{\rm d}\Omega'\, =\, \phi^*d\Omega
\end{equation*}
is the pullback to ${\mathbb P}_H$ of the symplectic form $\Omega$
on ${\mathcal M}^s_H$.

\subsection{Parabolic triples and spectral data}\label{para-spec.-da}
Take a parabolic triple $(E_*\, ,\theta\, ,\sigma)\, \in\,
{\mathbb P}_H$. Its image $\psi((E_*\, ,\theta))\, \in\,
{\mathcal H}$ will be denoted by $h$, where $\psi$ is the Hitchin map
defined in \eqref{def.-hitchin-map}.

As we saw in Section \ref{sec.-spectral-data},
the parabolic Higgs bundle $(E_*\, ,\theta)$ gives a spectral
curve $Y_h$ and a torsionfree sheaf $\mathcal L$ of rank one on
$Y_h$. It was noted in \eqref{direct-image-l} that
$\gamma_*{\mathcal L} \, \cong\, E$, where $\gamma$ is the
projection of $Y_h$ to $X$. Now, since $\gamma$ is a finite map, the
natural homomorphism
\begin{equation*}
H^i(Y_h,\, {\mathcal L}) \, \longrightarrow\, 
H^i(X,\, \gamma_*{\mathcal L})\, =\, H^i(X,\, E)
\end{equation*}
is an isomorphism for all $i\geq 0$. Therefore,
${\mathbb P}H^0(Y_h,\, {\mathcal L})
\,\cong \, {\mathbb P}H^0(X, E)$. The point in
${\mathbb P}H^0(Y_h,\, {\mathcal L})$ corresponding the point
$\sigma\, \in\, {\mathbb P}H^0(X, E)$ will be denoted
by $\sigma'$. In particular, $\sigma'$ is a divisor on $Y_h$.

We will calculate the degree of the divisor defined by
$\sigma'$. The Hitchin map $\psi$ in \eqref{def.-hitchin-map}
is a complete integrable system for the symplectic structure
${\rm d}\Omega$ on ${\mathcal M}_H$ and the fiber of $\psi$ over
any $h'\in {\mathcal H}$ is the Jacobian of the corresponding
spectral curve $Y_{h'}$. Consequently, the genus of
$Y_{h'}$ coincides with $\dim {\mathcal M}_H/2\, =\, 4g-3+n$. Since
\begin{equation*}
\text{degree}({\mathcal L}) -(4g-3+n)+1 \, =\,
\chi({\mathcal L}) \, =\, \chi(E) \, =\, d+2(1-g)
\end{equation*}
we conclude that $\text{degree}({\mathcal L}) \, =\,
d+n+2(g-1)$. Hence the degree of the divisor defined
by the section $\sigma'$ of $\mathcal L$, which coincides
with the degree of $\mathcal L$, is $d+n+2(g-1)$, where
$n\, =\, \# S$.

Set $\delta\, =\, d+n+2(g-1)$. Let $\text{Hilb}^\delta(Z)$ 
denote the Hilbert scheme of zero dimensional subschemes
of the surface $Z$ (the total space of
$K_X\otimes {\mathcal O}_X(S)$) of length $\delta$.

The divisor of $\sigma'$ is a zero dimensional subscheme of
$Y_h$ of length $\delta$. Taking the image of $\sigma'$ by the
inclusion map of the spectral curve $Y_h$ in $Z$, we get
a zero dimensional subscheme of $Z$ of length $\delta$.
Therefore, we have an element of
$\text{Hilb}^\delta(Z)$. Let
\begin{equation}\label{math-Hilb.-scheme}
f\, :\, {\mathbb P}_H\, \longrightarrow\, \text{Hilb}^\delta(Z)
\end{equation}
be the map that sends any parabolic triple to the zero
dimensional subscheme of $Z$ defined by the divisor of
the corresponding section on the spectral curve for the underlying
parabolic Higgs bundle. In other words, $f$ sends
any $(E_*\, ,\theta\, ,\sigma)\, \in\, {\mathbb P}_H$ to the
element of $\text{Hilb}^\delta(Z)$ defined by the divisor of
$\sigma'$ on $Y_h$.

\subsection{Forms on moduli of triples}\label{sec.-main-result}

Using the map $f$ defined in \eqref{math-Hilb.-scheme} we will
construct an algebraic one--form on ${\mathbb P}_H$, and for that
we will first show the existence of a natural meromorphic one--form
on $\text{Hilb}^\delta(Z)$.

We start by observing that
the variety $Z$ has a natural meromorphic one--form with pole, of
order at most one, along the divisor $\gamma^{-1}(S)$,
where $\gamma$, as before is the projection of $Z$ to $X$.
The $\gamma^*{\mathcal O}_X(S)$ valued one--form sends any tangent
vector $v\, \in\, T_z Z$ to $z\otimes{\rm d}\gamma(z)(v)$, where
${\rm d}\gamma(z)\, :\, T_z Z\, \longrightarrow\, T_{\gamma(z)}X$ is
the differential of $\gamma$ at $z$.
Note that since $z$ is an element of the
fiber $(K_X\otimes{\mathcal O}_X(S))_{\gamma (z)}$, the tensor
product $z\otimes d\gamma(z)(v)$ gives an element of the fiber
$({\mathcal O}_X(S))_{\gamma (z)}$ after contracting
$(K_X)_{\gamma (z)}$ with $T_{\gamma (z)}X$. Since $S$ is an
effective reduced divisor, a $\gamma^*{\mathcal O}_X(S)$ valued
one--form on $Z$ is a meromorphic one--form on $Z$ with a pole
along $\gamma^{-1}(S)$ of order at most one.

\textit{The meromorphic one--form on $Z$ defined
above will be denoted by $\Omega_Z$.} Using $\Omega_Z$, a meromorphic
one--form will be constructed on Hilbert
scheme $\text{Hilb}^k(Z)$ of zero dimensional subschemes of $Z$
of length $k$, where $k\, \geq \,1$.

Consider the Zariski open dense subset $U_k\, \subset\,
\text{Hilb}^k(Z)$
consisting of \textit{distinct} $k$ points of $Z$. Let
\begin{equation*}
\underline{z}\, = \, \{z_1, z_2, \cdots ,z_k\}\, \in\,
\text{Hilb}^k(Z)
\end{equation*}
be a point of $U_k$, that is, all $z_i$ are distinct. Then we have
\begin{equation*}
T_{\underline{z}} \text{Hilb}^k(Z)
\, =\, \bigoplus_{i=1}^k T_{z_i} Z\, .
\end{equation*}
Therefore, the meromorphic one--form $\Omega_Z$ on $Z$ defines a 
meromorphic one--form on the Zariski open subset $U_k$
of $\text{Hilb}^k(Z)$.
In other words, this form sends any tangent vector
\begin{equation*}
\{v_1, v_2, \cdots ,v_k\}\, \in\, T_{\underline{z}}
\text{Hilb}^k(Z)
\end{equation*}
where $v_i\in T_{z_i} Z$, to
\begin{equation*}
\sum_{i=1}^k {\Omega}_Z (z_i)(v_i)
\end{equation*}
whenever the sum makes sense. Evidently, this one--form is regular
on the complement of the divisor on
$U_k$ defined by all points $\{z_1, z_2, \cdots ,z_k\}$ such that
\begin{equation*}
\{\gamma(z_1), \gamma(z_2), \cdots ,\gamma(z_k)\}\cap S\, \not=\,
\emptyset
\end{equation*}
as the above sum makes sense on the complement. More precisely, the
pole of the one--form on $U_k$ is over this divisor, and the order
of the pole is one. Since $U_k$ is a Zariski open dense subset of
$\text{Hilb}^k(Z)$, the meromorphic one--form on $U_k$ defines a
meromorphic one--form on $\text{Hilb}^k(Z)$.

\textit{The meromorphic one--form on ${\rm Hilb}^k(Z)$
defined above will be denoted by $\Omega_k$.} Note that
$\text{Hilb}^1(Z) = Z$ and $\Omega_1 \, =\, \Omega_Z$.

Consider the meromorphic one--form $f^*\Omega_\delta$  on
${\mathbb P}_H$, where the map
$f$ is defined in \eqref{math-Hilb.-scheme}.
It was noted in Section \ref{sec.-spectral-data} that for any
spectral curve $Y_h$ and any $s\in S$, we have $\gamma^{-1}(s)\cap
Y_h\, =\, \{0\}$ (that is, the spectral curve is
totally ramified over $s$ and
passes through $0$). From this it follows immediately that
$f^*\Omega_\delta$ is a holomorphic one--form on ${\mathbb P}_H$.
Indeed, for the origin
$0\in (K_X\otimes{\mathcal O}_X(S))_x$, where $x\in X$,
the form $\Omega_Z$ vanishes at $0$. Therefore, $f^*\Omega_\delta$
is a holomorphic one--form on ${\mathbb P}_H$.

Recall the one--form $\Omega'\, =\, \phi^*\Omega$ on
${\mathbb P}_H$ constructed in \eqref{triple-one-form}.

\begin{theorem}\label{main-theroem}
The one--form $f^*\Omega_\delta$ on ${\mathbb P}_H$ coincides
with the one--form $\Omega'$. In particular,
${\rm d}f^*\Omega_\delta$ coincides with the pullback
$\phi^*{\rm d}\Omega$ of the symplectic form ${\rm d}\Omega$
on ${\mathcal M}^s_H$.
\end{theorem}

This theorem will be proved in the next section.

\section{Identification of one--forms}

We start with the following lemma.

\begin{lemma}\label{lemma1}
There is a holomorphic one--form $\omega$ on ${\mathcal M}^s_H$
such that $\phi^*\omega$ coincides with $f^*\Omega_\delta$,
where $\phi$ is the projection defined in \eqref{moduli-triple}.
\end{lemma}

\begin{proof}
Recall that ${\mathcal P}_H$ is a projective bundle over
${\mathcal M}^s_H$ with $\phi$ being the projection map. Since
there is no nonzero holomorphic one--form on a projective space,
the pullback of $f^*\Omega_\delta$ by the inclusion map of a
fiber of $\phi$ must vanish identically.

Let $T^\phi\, \subset \, T{\mathcal P}_H$ be the relative tangent
bundle. In other words, $T^\phi$ is the kernel of the differential
${\rm d}\phi \, :\, T{\mathcal P}_H\, \longrightarrow\,
\phi^*T{\mathcal M}^s_H$. Since the evaluation of
$f^*\Omega_\delta$ on $T^\phi$ vanishes, there is a homomorphism
\begin{equation*}
\omega'\, :\, T{\mathcal P}_H/T^\phi\, \longrightarrow\,
{\mathcal O}_{{\mathcal P}_H}
\end{equation*}
such that $f^*\Omega_\delta$ coincides with the composition of
the natural projection of $T{\mathcal P}_H$ to $T{\mathcal P}_H/T^\phi$
with the homomorphism $\omega'$; here ${\mathcal O}_{{\mathcal P}_H}$
is the structure sheaf of ${\mathcal P}_H$, or equivalently, the
sheaf defined by the trivial line bundle.

For any point $\zeta \, \in \, {\mathcal M}^s_H$, the restriction
of $T{\mathcal P}_H/T^\phi$ to the fiber $\phi^{-1}(\zeta)$ is a
trivial vector bundle. In fact, $T{\mathcal P}_H/T^\phi$
is identified with the pullback $\phi^*T{\mathcal M}^s_H$.
Since $\phi^{-1}(\zeta)$ is a compact and connected,
the homomorphism $\omega'(p)$ is independent of $p\, \in\,
\phi^{-1}(\zeta)$ (with $\zeta$ fixed). In other words, there
is a holomorphic one--form $\omega$ on ${\mathcal M}^s_H$ such that
$\omega'$ is the pullback of $\omega$. This completes the proof
of the lemma.
\end{proof}

Recall the one--form $\Omega$ on ${\mathcal M}^s_H$ constructed
in Section \ref{subsec.-symp.-str.}. Since $\Omega'\, =\,
\phi^*\Omega$ (see \eqref{triple-one-form}),
in view of Lemma \ref{lemma1}, to prove Theorem \ref{main-theroem}
it suffices to establish that the
two one--forms $\omega$ and $\Omega$ on ${\mathcal M}^s_H$
coincide.

Consider the Hitchin map $\psi$ defined in \eqref{def.-hitchin-map}.
We want to show that there is a holomorphic one--form $\Omega_H$
on the Hitchin space $\mathcal H$ such that
\begin{equation}\label{form-identity}
\Omega - \omega \, =\, \psi^*\Omega_H \, .
\end{equation}
To prove the existence of such a form $\Omega_H$, take a point
$h\,\in\, {\mathcal H}$ such that the corresponding spectral curve
$Y_h$ is smooth. We recall that there is a nonempty Zariski open
subset $U$ of $\mathcal H$ such that for any point $h'\,\in \, U$
the corresponding spectral curve $Y_{h'}$ is smooth.

The fiber $\psi^{-1}(h)$ is identified with the Picard variety
$J_h \,:= \, \text{Pic}^{d+2(1-g)}(Y_h)$ of degree
$d+2(1-g)$ line bundles on $Y_h$. (The degree $d+2(1-g)$ was
computed in the proof of Lemma \ref{lemma-vanish}.) Let
$j_h\, :\, J_h\, \longrightarrow\, {\mathcal M}^s_H$ be the
inclusion map of the fiber of $\psi$. From the
constructions of $\Omega$ and $\omega$ it follows that
that $j^*_h\Omega \, =\, j^*_h\omega$.

Therefore, exactly as in the proof of Lemma \ref{lemma1},
we conclude that for any point $z\in \psi^{-1}(h)$,
the homomorphism
\begin{equation*}
(\Omega - \omega)(z)\, :\, T_z{\mathcal M}^s_H \, \longrightarrow
\, {\mathbb C}
\end{equation*}
factors through the projection ${\rm d}\psi (z) \, :\,
T_z{\mathcal M}^s_H \, \longrightarrow\, T_{\psi (z)}{\mathcal H}$
defined the differential of $\psi$ at the point $h$. Consequently,
there is a holomorphic one--form $\Omega_H$ on $U$ such that
$\Omega - \omega \, =\, \psi^*\Omega_H$ on $\psi^{-1}(U)$, where
$U$, as before, is the open subset of $\mathcal H$
defined by the points corresponding
to smooth spectral curves. Since $\Omega - \omega$ on
$\psi^{-1}(U)$ extends to ${\mathcal M}^s_H$, $U$ is a Zariski
open dense subset of $\mathcal H$ and the map $\psi$ is a
submersion everywhere, it follows immediately
that $\Omega_H$ extends to $\mathcal H$ and the equality
in \eqref{form-identity} is valid on ${\mathcal M}^s_H$.

\begin{lemma}\label{lemma2}
The one--form $\Omega_H$ on $\mathcal H$ vanishes identically.
\end{lemma}

\begin{proof}
Recall that ${\mathcal H} = H^0(X, K_X)\oplus H^0(X, K^{\otimes
2}_X\otimes{\mathcal O}_X(S))$. Set ${\mathcal H}' \, :=\,
{\mathcal H}\setminus\{0\}$, the nonzero vectors. For any
nonzero complex number $c$, consider the automorphism of
${\mathcal H}'$ that sends any point $(\alpha\, ,\beta)$
to $(c\alpha\, ,c^2\beta)$, where $\alpha \in
H^0(X, K_X)$ and
\begin{equation*}
\beta \,\in\, H^0(X, K^{\otimes 2}_X\otimes{\mathcal O}_X(S))\, .
\end{equation*}
So we have a free
action of ${\mathbb C}^*$ on ${\mathcal H}'$ defined this way.
The quotient space
\begin{equation*}
Q\, :=\, {\mathcal H}'/{\mathbb C}^*
\end{equation*}
is a weighted projective space. Let
\begin{equation}\label{wt.-proj.-sp.}
\rho\, :\, {\mathcal H}' \, \longrightarrow Q
\end{equation}
be the quotient map.

In Section \ref{sec.-spectral-data}, given a point of $\mathcal H$
we constructed a spectral curve, which is a divisor on $Z$, the
total space of $K_X\otimes{\mathcal O}_X(S)$. We want to describe
the above action of ${\mathbb C}^*$ on ${\mathcal H}'$ in terms of
spectral curves. On $Z$ there is an action of ${\mathbb C}^*$
defined by the condition that the action of any $c\,\in\,
{\mathbb C}^*$
sends a point $z$ to $cz$, where the scalar multiplication
is defined by the vector space structure of the fibers of the
line bundle $K_X\otimes{\mathcal O}_X(S)$. It is easy to see that
for any $h\,\in\, {\mathcal H}'$ and
$c\,\in\, {\mathbb C}^*$, the spectral
curve corresponding to the point $ch$ coincides with the image of the
spectral curve corresponding to the point $h$ by the automorphism
of $Z$ defined action of $c$.

On the other hand, the meromorphic one--form $\Omega_Z$ on $Z$
(constructed in Section \ref{sec.-main-result}) evidently has
the property
that it vanishes along the orbits of ${\mathbb C}^*$ on $Z$.
In other words, for the projection $\gamma$ of $Z$ to $X$, the
pullback of $\Omega_Z$ by the inclusion map of a fiber of
$\gamma$ vanishes identically. Furthermore, for any $c\in
{\mathbb C}^*$, if $T_c$ denotes the automorphism of $Z$ defined
by the multiplication by $c$, then $T^*_c\Omega_Z \, =\,
c\Omega_Z$. From these observations it follows immediately
that there is a one--form $\Omega_Q$ on the weighted projective
space $Q$ such that
$\rho^*\Omega_Q \, =\, \Omega_H$ on ${\mathcal H}'$, where $\rho$
is the projection in \eqref{wt.-proj.-sp.}.

A weighted projective space does not admit any nonzero
holomorphic one--form. Hence we have $\Omega_Q\, =\, 0$.
Since $\rho^*\Omega_Q \, =\, \Omega_H$, it follows immediately
that $\Omega_H\, =\, 0$, and the proof of the lemma is complete.
\end{proof}

The above results clearly combine together in imply Theorem
\ref{main-theroem}.

\noindent {\it Proof of Theorem \ref{main-theroem}.} Lemma
\ref{lemma2} \eqref{form-identity} together imply that
$\Omega\, =\, \omega$ on ${\mathcal M}^s_H$. So,
\begin{equation*}
\Omega'\, :=\, \phi^*\Omega\, =\, \phi^*\omega\, .
\end{equation*}
Now Lemma \ref{lemma1}, which says that $\phi^*\omega\,
=\, f^*\Omega_\delta$, completes the proof of the theorem.


\end{document}